\documentclass[12pt]{article}

\usepackage{amssymb,amsmath,amsthm,latexsym}

\newtheorem{theorem}{Theorem}

\newtheorem{remark}[theorem]{Remark}
\newenvironment{rem}
 {\begin{remark}\normalfont}
 {\end{remark}}

\evensidemargin=\oddsidemargin

\begin{document}

$\;$

\vspace{.8in}

{\large{
\centerline{Constructing regular self-similar solutions to the 3D Navier-Stokes equations}
\centerline{originating at singular and arbitrary large initial data}}}

\vspace{.4in}

\centerline{
Zoran Gruji\'c\footnote{Department of Mathematics, University
of Virginia; zg7c@virginia.edu} }

\vspace{.4in}

\begin{abstract}
Global-in-time smooth self-similar solutions to the 3D Navier-Stokes equations
are constructed emanating from homogeneous of degree $-1$ \emph{arbitrary large}
initial data belonging only to the closure of the test functions in
$L^2_{loc,unif}$.

\end{abstract}



\newpage

\section{Introduction}

A very simplistic description of the existing global-in-time results for the 3D Navier-Stokes 
equations is as follows: weak solutions are constructed from arbitrary large (mostly) $L^2$-type 
initial data, while mild and strong solutions are generated by small initial data in a variety of 
functional spaces. It is not known if the weak solutions are regular and there are only partial
regularity results -- Caffarelli, Kohn and Nirenberg \cite{CKN} proved that in a class of weak solutions
satisfying \emph{localized} energy inequality, the so called \emph{suitable weak solutions}, for every time 
$T>0$ the one-dimensional Hausdorff measure of the (possible) singular set in $\Omega \times (0,T)$ is 0 
($\Omega$ denotes the spatial domain). The standard constructions of weak solutions including the original
Leray's construction \cite{L} and the construction of the suitable weak solutions \cite{CKN} start 
from the inital data belonging to $L^2$ on the whole spatial domain. Recently, Lemari\'e-Rieusset 
\cite{L-R1}
(henceforth, the references will be to the book \cite{L-R2})
constructed \emph{localized} weak solutions started from initial data in the closure of 
the test functions in $L^2_{loc, unif}$. 

Scaling-invariant solutions of any physical model are of special interest -- heuristically (and
some times rigorously), they capture the qualitative behavior of the model. It turns out that
there is a unique scaling invariance for the Navier-Stokes system: if $u$ is a solution on
$\mathbb{R}^3 \times (0,\infty)$ with initial data $u_0$
then $v^\lambda(x,t) = \lambda u(\lambda x, \lambda^2 t)$ is a solution corresponding to
the initial data $v^\lambda_0(x) = \lambda u_0(\lambda x)$ for any $\lambda > 0$. If the original
initial data $u_0$ are homogeneous of degree $-1$, i.e. if $u_0(\lambda x)=\lambda^{-1}u_0(x)$ for
all $x$ in $\mathbb{R}^3$ and all $\lambda >0$, and the problem is considered in a uniqueness
class for the Navier-Stokes equations, it follows readily that 
$u(\lambda x, \lambda^2 t)=\lambda^{-1} u(x,t)$ for all $x$ in $\mathbb{R}^3$, all $t>0$
and all $\lambda>0$, i.e.
$u$ is scaling-invariant or self-similar. It is worth noting that initial data generating (nontrivial)
self-similar solutions are necessarily singular -- the homogeneity of degree $-1$ implies 
the existence of a singularity at the origin of the order of $\frac{1}{|x|}$. This prevents
obtaining self-similar solutions in the `usual' spaces, e.g., the Lebesgue spaces $L^p$ for $p \ge 3$ 
or the Sobolev spaces $H^s$ for $s \ge 1/2$. In fact, the $\frac{1}{|x|}$-type singularity
is exactly a borderline singularity for the borderline spaces $L^3$ and $H^{1/2}$. The first
rigorous construction of self-similar solutions is due to Giga and Miyakawa \cite{GM}. They worked
in the Morrey-type spaces of measures and obtained self-similar solutions for the vorticity
equation starting from small initial data. More recently, Cannone \cite{C}
developed a systematic approach
to constructing self-similar solutions in a variety of spaces (including
the homogeneous Besov spaces) originating at small initial data. (See also
\cite{CK, CMP, CP, P1, P2}.)

It is well-known there are no nontrivial self-similar solutions satisfying the Leray's energy
inequality (see, e.g., \cite{GM}). Consequently, the standard constructions of weak solutions are 
not suitable for generating self-similar solutions. However, the weak solutions constructed
in \cite{L-R2} satisfy only the \emph{localized} energy inequality and hence are not \emph{a priori} 
incompatible with (nontrivial) self-similarity. 

In this paper we obtain existence of \emph{regular} self-similar 
solutions emanating from \emph{arbitrary large} initial data (homogeneous of degree $-1$) belonging
only to $L^2_{loc, unif}$. 
To the best of the author's knowledge, this is the first time global-in-time
existence of smooth solutions to the 3D Navier-Stokes equations is established 
without \emph{any} smallness condition on the initial data.

The proof starts off by
showing that a modified local-in-time part of the construction of weak solutions in 
\cite{L-R2} yields
``partially self-similar'' solutions on $(0,T)$.
Partial self-similarity is then used to infer spatial
regularity on $(0,T)$ via Caffarelli-Kohn-Nirenberg estimate on the size of a singular 
set and local in space-time Serrin's $L^p_x L^q_t$-criterion. 
Spatial continuity in turn leads to full self-similarity on $(0,T)$.
This interplay between self-similarity
and regularity is exploited further until fully self-similar, $C^\infty_x C^\infty_t$-solutions on
$\mathbb{R}^3 \times (0,\infty)$ are constructed.

\section{Preliminaries}

We consider the 3D Navier-Stokes equations (NSE) with the unit viscosity and the zero
external force,
\begin{align}\label{3d}
 &u_t - \triangle u + (u \cdot \nabla) u + \nabla p = 0\\
 &\nabla \cdot u = 0\notag
\end{align}
where $u$ denotes the velocity of the fluid and $p$ the pressure. The spatial domain 
will be $\mathbb{R}^3$. 

In what follows, a projected form of the equation will be useful:
\begin{align}\label{3dp}
 &u_t - \triangle u + P \, \nabla \cdot (u \otimes u) = 0\\
 &\nabla \cdot u = 0\notag
\end{align}
($P$ denotes the Leray projector). Let us mention here that if a distributional solution $u$
of (\ref{3dp}) on $\mathbb{R}^3 \times (0,T)$ is \emph{uniformly} locally square-integrable, then
there exists $p$ in $\mathcal{D}' \left(\mathbb{R}^3 \times (0,T)\right)$ such that
\[
 P \, \nabla \cdot (u \otimes u) = \nabla \cdot (u \otimes u) + \nabla p\
\]
in the sense of distributions, i.e. $(u,p)$ is a distributional solution of (\ref{3d}).

Lemari\'e-Rieusset defined a localized version of weak solutions, the so called
\emph{local Leray solutions}. The definition is as follows \cite{L-R2}.

Let $u_0 \in L^2_{loc, unif}$ with $\nabla \cdot u_0 = 0$.
A \emph{local Leray solution} $u$ corresponding to the initial data $u_0$ is a locally square-integrable 
distributional solution of (\ref{3d}) on $\mathbb{R}^3 \times (0,T)$ with
the following properties:

i) $u \in \cap_{t<T} L^\infty \left((0,t), L^2_{unif,loc}\right)$

ii) $\sup_{x_0 \in \mathbb{R}^3} \iint_{0<s<t, |x-x_0| \le 1} |\nabla \otimes u|^2 \; dx \; ds < \infty$
for all $t<T$

iii) $\lim_{t \to 0^+} \int_K |u-u_0|^2 \; dx = 0$ for any compact subset $K$ of $\mathbb{R}^3$

iv) $u$ is a \emph{suitable weak solution} in the sense of Caffarelli, Kohn and Nirenberg -- in
particular, $u$ satisfies the \emph{localized} energy inequality.

The following result can be found in \cite{L-R2} -- the idea of the proof is to use the localized
energy inequality as the main tool in constructing the solutions originating at the
$L^2_{loc,unif}$-initial data. Let $E$ denote the closure of the test functions in
the norm of $L^2_{loc,unif}$.

\begin{theorem}\label{l-r} \cite{L-R2} \ 
i) Let $u_0 \in L^2_{loc,unif}$ be divergence free. Then there exists a local Leray
solution $u$ on $\mathbb{R}^3 \times (0,T)$ for some $T>0$.

ii) If in addition $u_0 \in E$, then $u \in \cap_{t<T} L^\infty \left((0,t), E\right)$ and
$\lim_{t \to 0+} \|u-u_0\|_{L^2_{loc,unif}} = 0$.
\end{theorem}

The construction is based on deriving some (uniform in $\epsilon$)
local energy-type estimates for solutions
of the family of mollified Navier-Stokes equations
\begin{align}\label{3deps}
 &(u_\epsilon)_t - \triangle u_\epsilon + P \; \nabla \cdot \left((u_\epsilon * \rho_\epsilon)
  \otimes u_\epsilon\right) = 0\\
 &\nabla \cdot u_\epsilon = 0\notag
\end{align}
supplemented by $u_\epsilon(0, \cdot) = u_0$ for all $\epsilon > 0$.
A family of (standard) mollifiers
$\left\{\rho_\epsilon\right\}_{\epsilon>0}$
is defined as follows. Let $\rho$ be a non-negative test function normalized such that
$\int \rho \, dx = 1$. Then $\rho_\epsilon$ is given by 
$\rho_\epsilon(x)=\frac{1}{\epsilon^3} \rho(\frac{x}{\epsilon})$ for all $x$ in $\mathbb{R}^3$.
We will come back to this family of smooth approximating solutions in the following section.

\begin{rem}
It was shown in \cite{L-R2} that local-in-time local Leray solutions can be extended to
global-in-time local Leray solutions; however, we will not need that result here.
\end{rem}

Next, we recall a definition of (space-time) singular/regular points of a weak solution $u$ 
in the sense of \cite{CKN}.
A point $(x,t)$ in $\mathbb{R}^3 \times (0,T)$ is \emph{singular} if 
$\|u\|_{L^\infty(D)} = \infty$ for any space-time neighborhood $D$ of $(x,t)$. A point is
\emph{regular} if it is not singular, i.e. if there exists a space-time neighborhood $D$
such that $\|u\|_{L^\infty(D)} < \infty$.

The main partial regularity result is given by the following theorem.

\begin{theorem}\label{ckn} \cite{CKN} \
Let $u$ be a suitable weak solution on $\mathbb{R}^3 \times (0,T)$. Then the one-dimensional
Hausdorff measure of the possible singular set in $\mathbb{R}^3 \times (0,T)$ is 0. In
particular, the singular set can not contain a smooth space-time curve.
\end{theorem}

\begin{rem}
The original statement is somewhat stronger -- it states that the one-dimensional parabolic
Hausdorff measure of the singular set is 0.
\end{rem}

At first, it may be puzzling that local boundness is in this context identified as ``regularity''.
However, the following result due to Serrin \cite{S} implies smoothness in $x$ at any regular
point $(x,t)$.

\begin{theorem}\label{s} \cite{S} \
Let $u$ be a weak solution in some open space-time region $D$ such that
$u \in L^p_x L^q_t (D)$ for a pair of exponents $(p,q)$ satisfying the 
Foias-Prodi-Serrin condition
$\frac{3}{p}+\frac{2}{q}<1$. Then $u$ is of class $C^\infty$ in $x$, and each derivative
is bounded on compact subregions of $D$.
\end{theorem}

\section{The construction}

Since our construction of global-in-time regular self-similar solutions will originate 
at the construction of local-in-time local Leray solutions, the first thing to
check is whether the requirements on the initial data in Theorem \ref{l-r} are compatible
with the desired self-similarity of solutions. As already pointed out in the introduction, 
any nontrivial self-similar
solution must emanate from nontrivial initial data homogeneous of degree $-1$ and hence possessing
a $\frac{1}{|x|}$-type singularity at the origin. This is not in conflict with being
in $L^2_{loc,unif}$, and decay at infinity, i.e. a possibilty of an approximation with 
the test functions, is also compatible with the homogeneity of degree $-1$. A simple
example of a family of homogeneous of degree $-1$ initial data 
satisfying the conditions of Theorem \ref{l-r}
is the following:
\[
 u_0^\alpha(x)=\alpha\left(\frac{x_2-x_3}{|x|^2}, \frac{x_3-x_1}{|x|^2}, \frac{x_1-x_2}{|x|^2}\right)
\]
for $\alpha$ in $\mathbb{R}$. Notice that this family generates \emph{arbitrary large}
initial data in $L^2_{loc,unif}$.

The second imminent problem in constructing any type of invariant weak solutions is that
the weak solutions are not known to be unique, so although the inital data may be compatible
with the invariance, it is \emph{a priori} not clear if it is possible to construct 
any invariant solutions at all. One way out is to show that a particular construction of
weak solutions will preserve the invariance -- this was shown by Brandolese in the case
of a special rotational invariance in the \cite{CKN}-construction of suitable weak solutions
\cite{B1} (as pointed out in \cite{B2}, the same argument works for any rotational invariance,
and also for most of the standard constructions of weak solutions). 
The trick is to choose a test function in the definition of the smooth approximating 
solutions to be compatible with the invariance -- in this case any radial test function 
will do. The consequence is the invariance of the approximating solutions which is then
preserved in the limit. The same trick does not seem to work for the self-similarity. It turns
out that in order to construct self-similar approximating solutions, the test function 
would have to satisfy a certain degree of homogeneity which contradicts the very nature
of a test function, i.e. having a compact support. What will happen here is that although
the approximating solutions will \emph{not} be self-similar, just enough self-similarity
will be recovered in the limit to push the construction through.

\begin{theorem}\label{thm}
Let $u_0 \in E$ be divergence free and homogeneous of degree $-1$. Then there exists a
self-similar solution $u$ of (\ref{3dp}) on $\mathbb{R}^3 \times (0,\infty)$ corresponding to the 
initial data $u_0$ with the following properties:

i) $\lim_{t \to 0+} \|u(t)-u_0\|_{L^2_{loc,unif}} = 0$

ii) $u \in C^\infty_x C^\infty_t$ and hence a classical solution on $\mathbb{R}^3 \times (0,\infty)$

iii) $\|u(t)\|_{L^\infty} = \|u(1)\|_\infty \, \frac{1}{\sqrt{t}} < \infty$
for all $t$ in $(0,\infty)$.
\end{theorem}

\begin{proof}
Recall the construction of local-in-time local Leray solutions (\ref{3deps}). It is clear 
that (\ref{3deps}) will be scaling invariant if and only if the smoothed-out velocity
\[
 w_\epsilon(x,t)=(u_\epsilon * \rho_\epsilon)(x,t)=\int u_\epsilon(x-y,t) \rho_\epsilon(y) \; dy
\]
exhibits the same scaling as the original velocity, i.e. if and only if
$w_\epsilon(\lambda x, \lambda^2 t)=\lambda^{-1} w_\epsilon(x,t)$ for all $\lambda>0$.
A simple computation yields
$w_\epsilon(\lambda x, \lambda^2 t)=\lambda^{-1} 
\int u_\epsilon(x-z,t) \lambda^3 \rho_\epsilon(\lambda z) \; dz$
which is equal to $\lambda^{-1} w_\epsilon(x,t)$ if and only if 
$\rho_\epsilon(\lambda z)=\lambda^{-3} \rho_\epsilon(z)$ for all $\lambda>0$, i.e. 
if and only if $\rho_\epsilon$ is homogeneous of degree $-3$. This is not possible since $\rho_\epsilon$
is just a rescaled test function and hence have a compact support. Consequently, an approximating
solution $u_\epsilon$ will not be self-similar.

However there is some symmetry. For $0 < \lambda \le 1, 0 < \epsilon \le \epsilon_0$, define
$\rho^\lambda_\epsilon(x)=\lambda^3 \rho_\epsilon(\lambda x)$ and 
$u^\lambda_\epsilon(x,t)=\lambda u_\epsilon(\lambda x, \lambda^2 t)$ for all $x \in \mathbb{R}^3$ and
all $t \in (0,T)$.
A straightforward calculation reveals
that $u_\epsilon$ solves (\ref{3deps}) if and only if $u^\lambda_\epsilon$ solves
the following  system:
\begin{align}\label{3depslam}
 &(u^\lambda_\epsilon)_t - \triangle u^\lambda_\epsilon + 
  P \; \nabla \cdot \left((u^\lambda_\epsilon * \rho^\lambda_\epsilon)
  \otimes u^\lambda_\epsilon\right) = 0\\
 &\nabla \cdot u^\lambda_\epsilon = 0\notag
\end{align}
supplemented by $u^\lambda_\epsilon(0, \cdot) = u_0$ (utilizing the homogeneity of $u_0$).
Since $\rho^\lambda_\epsilon \ne \rho_\epsilon$
for any $\lambda \ne 1$, the approximating solutions $u^\lambda_\epsilon$ (for $\lambda \ne 1$) and
$u_\epsilon$ will differ as well. What saves the day is a simple observation that
$\rho^\lambda_\epsilon=\rho_{\frac{\epsilon}{\lambda}}$ and hence
$u^\lambda_\epsilon=u_{\frac{\epsilon}{\lambda}}$. Thus, $\{u_\epsilon\}$ is simply a subfamily of
$\{u^\lambda_\epsilon\}$. However, we need to be somewhat careful here. Since a weak solution $u$
in Theorem \ref{l-r} is obtained through a sequence, say $\{u_{\epsilon_n}\}$,
it is not possible to infer simultaneous convergence
for \emph{continuum} many lambda (and in particular,
for $\lambda$ in $(0,1]$). What we can do is to ``preconditon'' the sequence for 
a dense set. Let $D_2$ denote the set of all dyadics in $(0,1)$.
Playing with indices in
$D_2$, we can construct (via a diagonalization procedure) a sequence
$\{u_{\epsilon_n}\}$
with a property that
$\{u_{\frac{\epsilon_n}{\lambda}}\}$
converges to a solution $u$ for all $\lambda$ in $D_2 \cup \{1\}$.

Recall that one of the convergences obtained in the proof of Theorem \ref{l-r}
is a strong convergence of $\varphi u_{\epsilon_n}$ to $\varphi u$ in
$L^p\left((0,T), L^2(\mathbb{R}^3\right)$ for any test function
$\varphi \in \mathcal{D} (\mathbb{R}^3 \times (0,T))$, and any $p < \infty$. 
Fix a ball $B$ centered at the origin. Then, for any $\lambda$ in $D_2 \cup \{1\}$,
$a.e. \, t$,
\begin{equation}\label{lim1}
 \int_B u^\lambda_{\epsilon_n} (x,t) g(x) \, dx =
 \int_B u_{\frac{\epsilon_n}{\lambda}} (x,t) g(x) \, dx \to \int_B u(x,t) g(x) \, dx
\end{equation}
for all $g$ in $L^2(B)$.
On the other hand, 
it is easily seen that
\begin{equation}\label{lim2}
 \int_B u_{\epsilon_n}^\lambda (x,t) g(x) \, dx = 
 \int_B \lambda u_{\epsilon_n}(\lambda x,\lambda^2 t) g(x) \, dx
 \to \int_B \lambda u(\lambda x,\lambda^2 t) g(x) \, dx
\end{equation}
for all $g$ in $L^2(B)$.
Uniqueness of a limit coupled with (\ref{lim1})
and (\ref{lim2}) yields that for any $\lambda$ in $D_2$, $a.e. \, t$,
\[
 \int_B \left[u(x,t)-\lambda u(\lambda x,\lambda^2 t)\right] g(x) \, dx = 0
\]
for all $g$ in $L^2(B)$.

Now, fix $\lambda$ in $D_2$ and consider
\[
 f(t) = \int_B \left[u(x,t)-\lambda u(\lambda x,\lambda^2 t)\right] g(x) \, dx
\]
for $t$ in $(0,T)$. Since $f$ is continuous (by the weak time continuity of $u$
with values in $L^2(B)$)
and vanishes $a.e.$, $f$ is identically zero on $(0,T)$.
Next, for a fixed $t$ in $(0,T)$ consider
\[
 g(\lambda) = \int_B \left[u(x,t)-\lambda u(\lambda x,\lambda^2 t)\right] g(x) \, dx
\]
for $\lambda$ in $(0,1)$. This function is also continuous and vanishes 
on a dense set ($D_2$) -- hence vanishes for all $\lambda$ in $(0,1)$.
Finally, since the ball $B$ was arbitrary,
\begin{equation}\label{ss1}
 u(x,t)=\lambda u(\lambda x, \lambda^2 t)
\end{equation}
for any $\lambda$ in $(0,1)$, any $t$ in $(0,T)$, and $a.e. \ x$ in $\mathbb{R}^3$.
\emph{A priori}, the $a.e. \ x$-set depends on both $\lambda$
and $t$. One can actually show more uniformity,
but this will suffice for our construction.

Notice that we can extend (\ref{ss1}) to any $\lambda > 0$ provided 
that we restrict $t$ to \break
$(0,\min\{T,\frac{T}{\lambda^2}\})$.

Next, we show that the singular set of $u$ on $\mathbb{R}^3 \times (0,T)$ is empty. In what follows,
it will be convenient to define a family of parabolic space-time cylinders centered at a space-time 
point $(x,t)$. For $r>0$, define $C_{x,t}(r)$ by
\[
 C_{x,t}(r)=\{(y,\tau) \in \mathbb{R}^3 \times (0,\infty) | \; |\tau-t|<\frac{1}{2}r^2,
 |(x,\tau)-(y,\tau)|<r\}.
\]
We argue by contradiction. Suppose there exists a singular point $(x,t)$ of $u$ in
$\mathbb{R}^3 \times (0,T)$ and consider a family of parabolic cylinders $\{C_{x,t}(r)\}$
for $0< r \le \min\{\sqrt{t}, \sqrt{T-t}\}=R_{t,T}$. Since $(x,t)$ is a singular point,
$\|u\|_{L^\infty(C_{x,t}(r))}=\infty$ for all $r$. 
Fix any $\lambda$ in $\left(0, \sqrt{\frac{T}{t+\frac{1}{2}R_{t,T}^2}}\right)$, 
and then fix any $r$, $0 < r \le R_{t,T}$. For this choice of
$\lambda$, $u(\lambda y, \lambda^2 \tau)=\lambda^{-1} u(y,\tau)$ for any
$\tau$ in $(t-\frac{1}{2}r^2, t+\frac{1}{2}r^2)$, $a.e. \ y$ in $\mathbb{R}^3$, and hence
$a.e \ (y,\tau)$ in $C_{x,t}(r)$. This implies that
\[
 \|u\|_{L^\infty(C_{\lambda x, \lambda^2 t}(\lambda r))}
 = \lambda^{-1} \|u\|_{L^\infty(C_{x,t}(r))} = \infty
\]
for any $r$, $0 < r \le R_{t,T}$. In other words, we constructed a family of parabolic neighborhoods
of a point $(\lambda x, \lambda^2 t)$ on which $u$ blows-up and so $(\lambda x, \lambda^2 t)$
is also a singular point. The argument can be repeated for any 
$\lambda$ in $\left(0, \sqrt{\frac{T}{t+\frac{1}{2}R_{t,T}^2}}\right)$ and hence we 
generated a smooth space-time curve of singular points passing through the point $(x,t)$
contained in $\mathbb{R}^3 \times (0,T)$. Since local Leray solutions are suitable,
this contradicts Theorem \ref{ckn} and thus the singular set must be empty.

Now, it is standard to conclude infinite spatial regularity. Fix any space-time point $(x,t)$.
Since $(x,t)$ is a regular point, there exists a (bounded) neighborhood $D$ of $(x,t)$ such that
$\|u\|_{L^\infty(D)} < \infty$. This implies that $u \in L^p_x L^q_t (D)$ for any
$1 \le p,q, \le \infty$, and in particular $u$ satisfies the condition of Theorem \ref{s}.
It follows that $u \in C^\infty_x (D)$ and that each derivative is bounded on compact 
subregions of $D$.

Let us note that the spatial continuity of $u$ now implies that given $\lambda>0$,
\begin{equation}\label{ss2}
 u(x,t)=\lambda u(\lambda x, \lambda^2 t)
\end{equation}
for all $t$ in $(0, \min\{T, \frac{T}{\lambda^2}\})$ and all $x \in \mathbb{R}^3$.

Infinite time regularity in general does not follow from the emptiness of the
singular set. However, a restricted self-similarity will help. More precisely,
we can write $\displaystyle{u(x,t)=\frac{1}{\sqrt{t}} U_T\left(\frac{x}{\sqrt{t}}\right)}$
where $\displaystyle{U_T(\cdot)=\sqrt{\frac{T}{2}} 
u\left(\sqrt{\frac{T}{2}} \ \cdot, \frac{T}{2}\right)}$ for all
$(x,t)$ in $\mathbb{R}^3 \times (0,T)$. Since we know that $u \in C^\infty_x$, and it is clear that
$u \in C^\infty_x$ if and only if $U_T \in C^\infty$, it follows that $U_T$ is infinitely smooth.
The infinite time regularity is now readily seen.

To summarize, $u$ is in $C^\infty_x C^\infty_t$ (with the corresponding infinitely smooth
pressure $p$)
and hence a classical solutions on
$\mathbb{R}^3 \times (0,T)$ which is self-similar on $(0,T)$, i.e.
for a given $(x,t)$ in $\mathbb{R}^3 \times (0,T)$,
\begin{equation}\label{ss3}
 u(\lambda x, \lambda^2 t) = \lambda^{-1} u(x,t)
\end{equation}
for all $\lambda$ in $\left(0, \sqrt{\frac{T}{t}}\right)$ (the corresponding pressure 
scales as $p(\lambda x, \lambda^2 t)=\lambda^{-2}p(x,t)$).

Finally, we extend $(u,p)$ to the whole space-time by self-similarity. It is clear
that the extension will inherit infinite space and time regularity and will 
solve the equations (\ref{3d})
in the classical sense (pointwise).

One can now study various norms of $u$ -- we consider the $L^\infty$-norm here.
Recall that
$u \in \cap_{t<T} L^\infty\left((0,t), E\right)$ which paired with the time continuity
yields $u(t)$ in $E$ for any $t$ in $(0,T)$. In particular,
$u(\frac{T}{2})$ is in $E$ and hence it decays to $0$ as
$|x| \to \infty$. Since we know $u(\frac{T}{2})$ is continuous, it follows
$\|u(\frac{T}{2})\|_{L^\infty}<\infty$. Self-similarity implies
\[
 \|u(t)\|_{L^\infty} = \sqrt{\frac{T}{2}} \|u\left(\frac{T}{2}\right)\|_{L^\infty} \, \frac{1}{\sqrt{t}}
\]
for all $t$ in $(0,\infty)$. This is informative at both ends. It gives a blow-up rate as $t \to 0^+$
compatible with the $\frac{1}{|x|}$-type singularity of the initial data, and also a
generic decay rate as $t \to \infty$.

The convergence to the initial data in $L^2_{loc,unif}$ follows from
Theorem \ref{l-r}.

\end{proof}

\begin{rem}
It is worth noting that most of the functional spaces utilized in the study of
the 3D NSE on the whole space
are embedded in $L^2_{loc,unif}$ -- $u$ in $L^2_{loc}$ gives the bilinear term
$\nabla \cdot (u \otimes u)$ at least a distributional meaning, and uniformity
corresponds to the (desirable) translational invariance of a norm.
\end{rem}

\begin{rem}
It is clear that the proof of Theorem \ref{thm} implies infinite space and time
smoothness of \emph{any} self-similar suitable weak solution (regardless of its
construction).
\end{rem}

\vspace{1.5in}

ACKNOWLEDGMENTS \ The author thanks Prof. Ciprian Foias for his interest 
and insightful comments.


\vspace{2in}

\begin{thebibliography}{9999}

\bibitem[B1]{B1}
L. Brandolese, \emph{Asymptotic behavior of the energy and pointwise estimates for 
solutions to the Navier-Stokes equations} (preprint)

\bibitem[B2]{B2}
L. Brandolese, \emph{private communication}

\bibitem[CKN]{CKN}
L. Caffarelli, R. Kohn and L. Nirenberg, \emph{Partial regularity of suitable weak
solutions of the Navier-Stokes equations},
Comm. Pure Appl. Math. \textbf{35} (1982), 771-837

\bibitem[C]{C}
M. Cannone, \emph{Ondelettes, Paraproduits et Navier-Stokes}, Diderot Editeur, Paris, 1995

\bibitem[CK]{CK}
M. Cannone and G. Karch, \emph{Smooth or singular solutions to the Navier-Stokes
system ?} (preprint)

\bibitem[CMP]{CMP}
M. Cannone, Y. Meyer and F. Planchon, \emph{Solutions auto-similaires des \'equations
de Navier-Stokes dans $\mathbb{R}^3$}, Expos\'e n. VIII, S\'eminaire X-EDP,
Ecole Polytechnique, 1994

\bibitem[CP]{CP}
M. Cannone and F. Planchon, \emph{Self-similar solutions of the Navier-Stokes equations
in $\mathbb{R}^3$}, 
Comm. PDE \textbf{21} (1996), 179-193

\bibitem[GM]{GM}
Y. Giga and T. Miyakawa, \emph{Navier-Stokes flows in $\mathbb{R}^3$ with measures as 
initial vorticity and the Morrey spaces},
Comm. PDE \textbf{14} (1989), 577-618

\bibitem[L-R1]{L-R1}
P.G. Lemari\'e - Rieusset, \emph{Solutions faibles d'\'energie infinie pour les
\'equations de Navier-Stokes dans $\mathbb{R}^3$},
C.R. Acad. Sci. Paris S\'er. I \textbf{328 \, (8)} (1999), 1133-1138

\bibitem[L-R2]{L-R2}
P.G. Lemari\'e - Rieusset, \emph{Recent developments in the Navier-Stokes problem},
Research Notes in Mathematics, Chapman $\&$ Hall / CRC, 2002

\bibitem[L]{L}
J. Leray, \emph{Sur le mouvement d'un liquide visqueaux emplissant l'espace},
Acta Math. \textbf{63} (1934), 193-248

\bibitem[P1]{P1}
F. Planchon, \emph{Convergence de solutions des \'equations de Navier-Stokes vers de
solutions auti-similaires},
Expos\'e n. III, S\'eminaire X-EDP, Ecole Polytechnique, 1996

\bibitem[P2]{P2}
F. Planchon, \emph{Asymptotic behavior of global solutions to the Navier-Stokes equations
in $\mathbb{R}^3$},
Rev. Mat. Iberoamericana \textbf{14 \, 1} (1998), 71-93

\bibitem[S]{S}
J. Serrin, \emph{On the interior regularity of weak solutions of the Navier-Stokes
equations},
Arch. Rat. Mech. Anal. \textbf{9} (1962), 187-195

\end{thebibliography}
\end{document}